\documentclass[11pt]{amsart}
\usepackage{ amsmath}
\usepackage{amssymb}
\usepackage{amsxtra}
\usepackage{enumerate}
\usepackage[mathscr]{eucal}
\usepackage[margin=2.9 cm]{geometry}

\def\R{\mathbb R}

\def\H{\mathbb H}
\def \hh{{\bf H}_{\mathbb H}}

\def \fh{{\bf H}_{\F}}
\def \h{{\bf H}_{\mathbb R}}
\def\Sp{\mathrm{Sp}}
\def\SU{\mathrm{SU}}

\def\SL{{\rm SL}}

\def\SL{{\rm SL}}
\def\C{\mathbb C}

\def\PSL{{\rm PSL}}

\def\F{\mathbb F}

\def \s{\mathbb S}
\def\c{{\rm SL}(2, C_n)}

\newtheorem{theorem}{Theorem}[section]
\newtheorem{lemma}[theorem]{Lemma}

\theoremstyle{definition}
\newtheorem{definition}{Definition}

\theoremstyle{remark}
\newtheorem{remark}[theorem]{Remark}
\numberwithin{equation}{section}
\theoremstyle{plain}

\newtheorem{cor}[theorem]{Corollary}

\newcommand{\secref}[1]{Section~\ref{#1}}
\newcommand{\thmref}[1]{Theorem~\ref{#1}}
\newcommand{\lemref}[1]{Lemma~\ref{#1}}

\newcommand{\corref}[1]{Corollary~\ref{#1}}
\newcommand{\eqnref}[1]{~{\textrm(\ref{#1})}}

\begin{document}
\title[Discreteness of  Hyperbolic Isometries by Test Maps]{Discreteness Of Hyperbolic Isometries by Test Maps }
\author[K. Gongopadhyay]{Krishnendu Gongopadhyay}
 \address{ Indian Institute of Science Education and Research (IISER) Mohali,
Knowledge City, Sector 81, SAS Nagar, Punjab 140306, India}
\email{krishnendu@iisermohali.ac.in, krishnendug@gmail.com}

\author[A. Mukherjee]{Abhishek Mukherjee}
\address{ Kalna College, Kalna,  Dist. Burdwan, West Bengal 713409}
\email{abhimukherjee.math10@gmail.com }
\author[D. Tiwari]{Devendra Tiwari}
\address{Bhaskaracharya Pratisthana, 56/14, Erandavane, Damle Path,
Off Law College Road, Pune 411004. }
\email{devendra9.dev@gmail.com}

\subjclass[2000]{Primary 20H10; Secondary 51M10, 20H25 }
\keywords{ J\o{}rgensen inequality, discreteness, hyperbolic space, Clifford matrices.}
\date{\today}
\thanks{Gongopadhyay acknowledges partial support from SERB MATRICS grant MTR/2017/000355 and DST  grant DST/INT/RUS/RSF/P-19.}
\thanks{Mukherjee acknowledges partial support from the grant DST/INT/RUS/RSF/P-19.}
\thanks{Tiwari acknowledges support from the ARSI Foundation.}

\begin{abstract}
Let $\F=\R$, $\C$ or the Hamilton's quaternions $\H$.
Let $\fh^n$ denote the $n$-dimensional $\F$-hyperbolic space. Let ${\rm U}(n,1; \F)$ be the linear group that acts by the isometries of $\fh^n$.  A subgroup $G$ of ${\rm U}(n,1; \F)$  is called \emph{Zariski dense} if it does not fix a point on $\fh^n \cup \partial \fh^n$ and neither it preserves a totally geodesic subspace of $\fh^n$. 
We prove that a Zariski dense subgroup $G$ of ${\rm U}(n,1; \F)$ is discrete if for every loxodromic element $g \in G$, the two generator subgroup $\langle f, g  \rangle$  is discrete, where $f \in {\rm U}(n,1; \F)$ is a test map not necessarily  from $G$. 
\end{abstract}

\maketitle

\section{Introduction}\label{intro} 

Let $\F=\R$, $\C$ or the Hamilton's quaternions $\H$.   Let $\fh^n$ be the $n$-dimensional hyperbolic space over $\F$. Let ${\rm U}(n,1; \F)$ the unitary group that acts on $\fh^n$ by isometries. For simplicity of notations, ${\rm U}(n,1; \R)$ will be considered as the identity component of the full isometry group.  Following standard notations, we denote ${\rm U}(n,1; \R)={\rm{PO}}(n,1)$,   ${\rm U}(n,1; \C)={\rm U}(n,1)$, ${\rm U}(n,1; \H)=\Sp(n,1)$. 

The J\o{}rgensen inequality is an important result on discreteness of subgroups in two and three dimensional real hyperbolic geometry.   It was developed by J\o{}rgensen and later generalized to arbitrary dimension by Martin \cite{martin} and Waterman \cite{wat} using different approaches.  Abikoff and Haas \cite{ah}  proved that a Zariski-dense subgroup $G$ of ${\rm{PO}}(n,1)$ is discrete if and only if every two-generator subgroup of $G$ is discrete, also see \cite{martin}, \cite{fn}, \cite{lw}, \cite{wlc}. Following this theme,   Chen, in \cite{chen}, has obtained a discreteness criterion that uses a fixed `test map' to check  discreteness of a  subgroup.  Chen proved that a Zariski-dense subgroup $G$ of ${\rm{PO}}(n,1)$ is discrete if for each $g$ in $G$, the group $\langle g, h \rangle$ is discrete, where $h$ is a fixed non-trivial element from ${\rm{PO}}(n,1)$, not necessarily from $G$, such that $h$ is either of infinite order but not an irrational rotation, or if having finite order, it does not pointwise fix the minimal sphere containing the limit set of $G$. 

Chen's work suggests that the discreteness is not completely an internal property of a subgroup $G$, and one may detect it by performing discreteness of the two-generator subgroups having a fixed generator that might also be an element in the complement of $G$.  Such a generator is called a `test map'. The action of $\SL(2, \C)$  on the Riemann sphere by the linear fractional transformations provides an identification of ${\rm{PO}}(3,1)$ with $\PSL(2, \C)$.  In \cite{yang1}, \cite{yang2} and \cite{cao},  refined versions of discreteness criteria in $\SL(2, \C)$ using test maps have been obtained.  A generalization of the complex linear fractional transformations are the quaternionic linear fractional transformations that can be identified with the group $\PSL(2, \H)$. Here $\SL(2, \H)$, the $2 \times 2$ quaternionic matrices with Dieudonn\'e determinant $1$, acts by the linear fractional transformations on the boundary of the $5$-dimensional hyperbolic space.  In \cite{gm}, also see \cite{kel},  some  J{\o}rgensen type inequalities for two generator subgroups of $\SL(2, \H)$ were obtained. In \cite{gm2}, these inequalities are used to prove that the discreteness of a Zariski-dense subgroup $G$ of $\SL(2, \H)$ is determined by the two generator subgroups $\langle f, g \rangle$, where $f$ is a certain test map from $\SL(2, \H)$ and $g$ is a loxodromic element of $G$.

 In this short note, we generalise the above results to ${\rm U}(n,1; \F)$  to show that the discreteness of a subgroup $G$ is determined by a test map and the loxodromic elements of $G$. We also provide some quantitative bounds for the test maps.

\medskip Recall that an element $f$ in ${\rm U}(n,1; \F)$ is called \emph{elliptic} if it has a fixed point on $\fh^n$, it is \emph{parabolic}, resp. \emph{loxodromic} (or hyperbolic) if it has exactly one, resp. two fixed points on $\partial \fh^n$ and no fixed point on $\fh^n$. An elliptic element is called \emph{regular} if it has a unique fixed point on $\fh^{n}$.   This  type of isometries exist in  all dimensions over $\F=\C, \H$. However,  regular elliptic isometries of $\h^n$ exist if and only if $n$ is even.  When $n$ odd, every elliptic isometry of $\h^n$ has at least two fixed points on the boundary $\partial \h^n$. By abuse of notation, an elliptic isometry of $\h^n$ will be called \emph{regular} if it has at most two boundary fixed points.   A subgroup $G$ of ${\rm U}(n,1; \F)$ is called \emph{Zariski-dense} or \emph{irreducible} if it does not have a global fixed point on $\overline \fh^{n}=\fh^{n} \cup \partial \fh^{n}$ and neither it preserves a proper totally geodesic subspace of $\fh^{n}$.

\subsection{Discreteness in ${\rm{PO}}(n,1)$} For ${\rm {PO}}(n,1)$ we  use the Clifford algebraic formalism that  was initiated by Ahlfors in \cite{ahl}, \cite{ahl2}. Waterman gave an alternative formulation of this approach in \cite{wat} and proved its equivalence to Ahlfors's formalism. In this approach the Clifford group $\SL(2, C_n)$, $n \geq 0$,  acts by the orientation-preserving isometries of $\h^{n+2}$, $n \geq 0$. The action is by the familiar looking linear fractional transformations. The group $\c$ consists of the $2 \times 2$ invertible matrices over Clifford numbers with  `Clifford determinant' one.  Waterman obtained J\o{}rgensen type inequalities for two-generator subgroups of  $\SL(2, C_n)$ in \cite{wat}. 

Cao and Waterman extended Waterman's inequalities using conjugacy invariants in \cite{cw}. Given an isometry $f$ of $\h^{n+2}$, one can associate `rotation angles' to it, and the rotation angles may be chosen to be elements of $(-\pi, \pi]$. The rotation angles are conjugacy invariants of an element and  one can further classify dynamical types of elements in $\c$  using the rotation angles and translation lengths, see  \cite{goku}. For a non-elliptic isometry $f$, let $\tau_f$ denotes the translation length of $f$ between the fixed points. $\tau_f=0$ if and only if $f$ is parabolic. The conjugacy invariant $\beta(f)$ used by Cao and Waterman can be defined as follows. 

 \begin{definition}\label{betaf}
Let $f$ be an element in $\SL(2, C_n)$. Let $\theta_1, \ldots, \theta_k \in (-\pi, \pi]$ be rotation angles of $f$ (counted with multiplicities). Let $\Theta=\max_{1 \leq i \leq k} |\theta_i|$.

 If $f$ is elliptic or parabolic, then $\beta(f)=4  \sin^2(\Theta/2)$.

If $f$ is loxodromic, then
$\beta(f)=4 \sinh^2(\tau_f/2)+4  \sin^2(\Theta/2)$.
\end{definition}
We apply the J\o{}rgensen type inequalities of Cao and Waterman to obtain discreteness criteria of a  Zariski-dense subgroup $G$ of $\SL(2, C_n)$ using test maps. We prove the following.
\begin{theorem}\label{thm1}
Let $G$ be a  Zariski-dense subgroup of $\c$. 
\begin{enumerate}
\item Let $f$ be a loxodromic element in $\c$,  not necessarily in $G$, such that $0< \beta(f)< 1$.  If  the two generator subgroup $\langle f,g\rangle$ is discrete for every loxodromic element $g$ in $G$,  then $G$ is  discrete.

\medskip \item Let $f$ be a non-elliptic isometry in $\c$, not necessarily in $G$, such that $$0< 2\cosh(\tau_f/2)\sqrt{\beta(f)} < 1.$$ If  the two generator subgroup $\langle f,g \rangle$ is discrete for every loxodromic element $g$ in  $G$,   then $G$ must be discrete.

\medskip \item Let $f$ be an elliptic element in $\c$, not necessarily in $G$, such that   $0< \beta(f) <4 \sin^2 ({\pi}/{10})$. If  the two generator subgroup $\langle f,g\rangle$ is  discrete for every loxodromic element $g$ in $G$,  then $G$ is  discrete.
\end{enumerate}
\end{theorem}

\medskip The following theorem also follows using similar methods as in the proof of the above theorem.
\begin{theorem}\label{thm2}
Let $G$ be a Zariski-dense subgroup of $\c$.
\begin{enumerate}
\item Let $f$ be a loxodromic element in $\c$,  not necessarily in $G$, such that $0< \beta(f) <1$.  If  the two generator subgroup $\langle f,gfg^{-1}\rangle$ is discrete for every loxodromic element $g$ in $G$,  then  $G$ is  discrete.
\item Let $f$ be a non-elliptic isometry in $\c$, not necessarily in $G$, such that $$0< \rho=2\cosh(\tau_f/2)\sqrt{\beta(f)} < 1.$$ If the two generator subgroup $\langle f,gfg^{-1} \rangle$ is discrete  for every loxodromic element $g$ in  $G$, then $G$ is discrete.

\item Let $f$ be a regular  elliptic element in  $\c$, not necessarily in $G$, such that \\ $0< \beta(f) <4 \sin^2 ({\pi}/{10})$.  If  the two generator subgroup $\langle f,gfg^{-1}\rangle$ is  discrete  for every loxodromic element $g$ in  $G$,  then  $G$ is  discrete.
\end{enumerate}
\end{theorem}

\subsection{Discreteness in ${\rm U}(n, 1; \F)$}  It is natural to ask for extending the above results to isometries of the complex and the quaternionic hyperbolic spaces.  Some discreteness criteria in $\SU(n,1)$ are available in the literature, eg. \cite{fh}, \cite{qjb}, \cite{qj}. However, not much attention has been given to $\Sp(n,1)$, partly because it lacks conjugacy invariants (unlike the  complex case) due to non-commutativity of the quaternions.  In the following we note a version of \thmref{thm1} in this set up.

A loxodromic element in $\Sp(n,1)$ is conjugate to a matrix of the form 
\begin{equation}\label{elpl}
f=diag(\lambda_1, \bar \lambda_1^{-1}, \lambda_3 \ldots, \lambda_{n+1}), 
\end{equation}
where $|\lambda_1|>1$, and $|\lambda_i|=1$ for $i=3, \ldots, n+1$.  Cao and Parker defined the following conjugacy invariant  in \cite{cp1}: 
$$\delta_{cp}(f)=\max\{|\lambda_i-1| \ : \ i=3, \ldots, n+1\}, $$ 
$$M_f=2 \delta_{cp}(f)+|\lambda_1-1|+|\bar \lambda_1^{-1}-1|.$$

An eigenvalue $\lambda$ of a matrix in $\Sp(n,1)$ is called negative-type or positive-type according as the Hermitian length of the corresponding eigenvector is negative or positive. An elliptic element in $\Sp(n,1)$  is conjugate  to a matrix of the form 
\begin{equation}\label{elp}
f=diag(\lambda_1, \ldots, \lambda_{n+1}), 
\end{equation}
where for all $i$, $|\lambda_i|=1$, and we choose the underlying Hermitian form so that $\lambda_1$ is a negative-type eigenvalue and all others are positive-type eigenvalues.  In \cite{gmt2}, we defined the following invariant,  cf. \cite{ct1}, 
\begin{equation} \label{eqq1} \delta(f)=\max\{|\lambda_i-1|+|\lambda_1-1| \ : \ i=2, \ldots, n+1\}. \end{equation} 
Clearly, $\delta(f)$ is an invariant of the conjugacy class of the elliptic element $f$.

Let $T_{s, \zeta}$ be a unipotent parabolic element in $\Sp(n,1)$. We shall call such element in 
$\Sp(n,1)$ or $\SU(n,1)$ as Heisenberg translation. 
We may assume (see \cite[p. 70]{chen}) that up to conjugacy, 
\begin{equation}\label{ht1} 
T_{s, \zeta}=\begin{pmatrix} 1 & 0 & 0 \\ s & 1 & {\zeta}^{\ast} \\ \zeta & 0 & I \end{pmatrix},
\end{equation}
where $Re(s)=\frac{1}{2} |\zeta|^2$. 

 \begin{theorem}\label{thmq} 
Let $G$ be Zariski dense in $\Sp(n, 1)$. 
\begin{enumerate}
\item Let $f \in \Sp(n, 1)$ be a loxodromic element such that $M_f < 1$. If $\langle f, g \rangle$ is discrete for every loxodromic element $g \in G$, then $G$ is discrete.

\item Let  $f \in \Sp(n, 1)$ be a  Heisenberg translation such that $|\zeta| < \frac{1}{2}$. If $\langle f,  g \rangle$ is discrete for every loxodromic element  $g$ in $G$, then $G$ is discrete. 

\item 
Let $f \in \Sp(n, 1)$ be a regular elliptic element such that  $\delta(f)<1$.  If $\langle f, g \rangle$ is discrete for every loxodromic element $g \in G$, then $G$ is discrete.

\end{enumerate} 
\end{theorem} 
As a by-product of the proof of the above theorem, we have the following result for subgroups in $\SU(n,1)$. A version of this result was obtained by Qin and Jiang in \cite{qj}.  
\begin{cor}\label{cord} 
Let $G$ be Zariski dense in $\SU(n, 1)$. 
\begin{enumerate}

\item Let $f \in \SU(n, 1)$ be a loxodromic element such that $M_f < 1$. If $\langle f, g \rangle$ is discrete for every loxodromic element $g \in G$, then $G$ is discrete.

\item Let  $f \in \SU(n, 1)$ be a  Heisenberg translation such that $|\zeta| < \frac{1}{2}$. If $\langle f,  g \rangle$ is discrete for every loxodromic element  $g$ in $G$, then $G$ is discrete. 
\item 
Let $f \in \SU(n, 1)$ be a regular elliptic element such that  $\delta(f)<1$.  If $\langle f, g \rangle$ is discrete for every loxodromic element $g \in G$, then $G$ is discrete.

\end{enumerate} 
\end{cor}

\medskip 
After discussing some background materials in \secref{prel}, we prove \thmref{thm1} and \thmref{thm2} in \secref{ms}. We prove \thmref{thmq} in \secref{qc}. 

\section{Preliminaries}\label{prel}
\subsection{Clifford Algebra}  The {\it Clifford algebra $C_n$}, $n \geq 0$,  is the real associative algebra which has been generated by $n$ symbols $ i_1,i_2,\cdots,i_n$ subject to the following relations:$$i_t i_s=-i_s i_t\;, \text{for}\; t \neq s \;\; \text{and}\; i_t^2=-1\;. $$
Let us define $i_0=1$ and then every element of $C_n$ can be expressed uniquely in the form $a= \sum a_I I$ , where the sum is over all products $I=i_{v_1}i_{v_2}\cdots i_{v_k}\;,$with $1\leq v_1 <v_2< \cdots <  v_k\leq n$ and $a_I \in \R$. Here the null product is permitted and identified with the real number $1$ . We equip $C_n$ with the Euclidean norm. Thus $C_0=\R$, $C_1=\C$, $C_2=\H$ etc. The following are involutions in $C_n$:
\medskip

$*$: In $a \in C_n$ as above, replace in each $I=i_{v_1}i_{v_2}\cdots i_{v_k}$ by $i_{v_k}\cdots i_{v_1}$. $a \mapsto a^{\ast}$ is an anti-automorphism.

$'$: Replace $i_k$ by $-i_k$ in $a$ to obtain $a'$.

The conjugate $\bar a$ of $a$ is now defined as: $\bar a=(a^{\ast})'=(a')^{\ast}$.

\medskip Let us identify  $\R^{n+1}$ with the $(n+1)-$dimensional subspace of $C_n$ formed by the Clifford numbers of the form
$$v=a_0 + a_1 i_1 +\ldots+a_n i_n.$$ These numbers are known as \emph{vectors}. The products of non-zero vectors form a multiplicative group,   denoted by $\Gamma_n$. For a vector $v$,  $v^{-1}=\bar v/|v|^2$.

\medskip A Clifford matrix of dimension $n$ is a $2 \times 2$ matrix $T=\begin{pmatrix} a & b \\ c & d \end{pmatrix}$ such that

(i) $a$, $b$, $c$, $d \in \Gamma_n -\{0\}$;

(ii) the Clifford determinant $\Delta(T)=  ad^{\ast}-bc^{\ast}=1$, and,

(iii) $ab^{\ast}$, $cd^{\ast}$, $c^{\ast} a$, $d^{\ast} b \in \R^{n+1}$.

\medskip The group of all Clifford matrices is denoted by $\SL(2, C_n)$.  In \cite{wat}, Waterman showed that $\SL(2, C_n)$ is same as the group of all invertible $2 \times 2$ matrices over $C_n$ with Clifford determinant 1.

The group $\c$ acts on $\s^{n+1}=\R^{n+1} \cup \{\infty\}$ by the action:
$$A: v \mapsto (a v+b)(c v + d)^{-1}. $$
This action extends by Poincar\'e extension to $\h^{n+2}$. The group $\c$ acts as the orientation-preserving isometry group of $\h^{n+2}$. 
For more details we refer to \cite{ahl}, \cite{ahl2}, \cite{wat}, \cite{cw}.

\subsection{Classification of elements of $\c$ :}\label{class}
We recall that, see \cite{wat},   a parabolic $f$ element in $\c$  is conjugate to
$$\begin{pmatrix} \lambda & \mu \\ 0 & {\lambda^*}^{-1} \end{pmatrix}, ~ |\lambda|=1, ~\mu \neq 0. $$
If $\lambda=1$, then $f$ is called a \emph{translation}.

Up to conjugacy in $\c$, a loxodromic element $f$  is  given by  $$f=\begin{pmatrix} \lambda & 0 \\ 0 & {\lambda^*}^{-1} \end{pmatrix},$$
where $\lambda \in \Gamma_n$, $|\lambda| \neq 1$. If $|\lambda|=1$, then it is a non-regular elliptic element.

Suppose $f$ is regular elliptic in $\c$, where $n$ is even. Note that $\SL(2, C_n)$ has a natural inclusion in $\SL(2, C_{n+1})$ as a closed subgroup. We shall consider the inclusion of $f$ in $\SL(2, C_{n+1})$,  and assume that $f$ fixes at least two points on the boundary $\partial \h^{n+3}$. Otherwise, we can choose two  fixed points of $f$ on $\partial \h^{n+2}$.   So, up to conjugacy in ${\rm SL}(2, C_{n+1})$, $f$ is of the form
$$\begin{pmatrix} \lambda & 0 \\ 0 & {\lambda^*}^{-1} \end{pmatrix}, ~|\lambda|=1.$$
The diagonal element $\lambda$  depends on the rotation angles of $f$, for details see \cite[Section 4]{wat}.
\subsection{Clifford Cross Ratio}
As in the complex analysis, Clifford cross ratios are defined similarly.
Let $z_1, z_2, z_3, z_4 \in \partial \h^{n+2}$ be any four distinct points. Let $z_1 \neq \infty$. The Clifford cross ratio of $(z_1, z_2, z_3, z_4)$ is given by
\begin{eqnarray*}
  [z_1, z_2, z_3, z_4]&=& (z_1-z_3)(z_1-z_2)^{-1}(z_2-z_4)(z_3-z_4)^{-1},\; \text{if}\; z_2, z_3, z_4 \neq \infty; \\
  &=& (z_1-z_3)(z_3-z_4)^{-1},\; \text{if}\; z_2=\infty; \\
  &=& (z_1-z_2)^{-1}(z_2-z_4),\; \text{if}\; z_3=\infty;\\
  &=& (z_1-z_3)(z-z_2)^{-1},\; \text{if}\; z_4=\infty.
\end{eqnarray*}

\medskip One can easily prove that for any $f=\begin{pmatrix}
          a & b \\
          c & d
        \end{pmatrix}\in \c,$ we have $$[fz_1,fz_2,fz_3,fz_4]={(cz_3+d)^*}^{-1}[z_1,z_2,z_3,z_4](cz_3+d)^*.$$

Thus $|[z_1,z_2,z_3,z_4]|$ and $\rm Re[z_1,z_2,z_3,z_4]$ are invariants of M\"{o}bius maps in $\c$. We have the following basic properties of cross ratios, see \cite{cw} for details.
\begin{enumerate}
  \item $[z_1,z_2,z_3,z_4]+[z_2,z_1,z_3,z_4]=1.$
  \item $[z_1,z_2,z_3,z_4][z_4,z_2,z_3,z_1]=1.$
  \item $|[z_1,z_2,z_3,z_4]|=|[z_2,z_1,z_4,z_3]|.$
  \item $|[z_1,z_2,z_3,z_4]|=|[z_3,z_4,z_1,z_2]|.$
\end{enumerate}

\subsection{Cao-Waterman J{\o}rgensen inequality}
We need call the following results which are important J{\o}rgensen type inequalities for two-generator subgroups of $\c$ when one of the generators is either elliptic or loxodromic.

\begin{theorem}\label{watt1} \cite{cw}
Let $ g= \begin{pmatrix}
          a & b \\
          c & d
        \end{pmatrix}\in \c$ be any element and $f  \in \c$ be a loxodromic element having two fixed points $u,v$ in $\partial \h^{n+2}$ satisfying that $\{gu,gv\}$ is not equal to $\{u,v\}$. If $\langle f,g \rangle$ generate a discrete subgroup in $\c$, then
$$\beta(f)\big (1+|[u,v,gu,gv]| \big) \geq 1.$$
\end{theorem}
\begin{theorem}\label{watt2} \cite{cw}
If $ g= \begin{pmatrix}
          a & b \\
          c & d
        \end{pmatrix} \in \c$ any element and $f  \in \c$ be an elliptic element such that $\langle f,g \rangle$ forms a non-elementary discrete subgroup in $\c,$ then we have
 $$\beta(f)\bigg (\frac{1}{4 \sin^2(\pi/10)}+|[u,v,gu,gv]| \bigg ) \geq 1,$$ where $u,v$ are any two boundary fixed points of $f$.
\end{theorem}
The J\o{}rgensen type inequality for non-elliptic isometries fixing the boundary point $\infty$ is given by the following.
\begin{theorem}\cite{cw} \label{watt3}
$f =\begin{pmatrix} \lambda & \mu\\ 0 & {\lambda^{*}}^{-1}\end{pmatrix}\in \c$ be a non-elliptic isometry that fixes the boundary point $\infty$. Let
Let $ g= \begin{pmatrix}
          a & b \\
          c & d
        \end{pmatrix}\in \c$ be any element in $\c$ such that  $0 < \rho=2\cosh(\tau_f/2)\sqrt{\beta(f)}< 1$,  and $fix(f) \cap fix(g) =\emptyset$. If $\langle f,g \rangle$ generate a discrete subgroup in $\c$, then
        $$|tr^2(fgf^{-1})[fg(\infty),fg^{-1}(\infty),g(\infty),g^{-1}(\infty)]| \geq  \frac{(1-\rho+\sqrt{(1-\rho)^2-4\beta(f)}}{2}\;.$$\\
Moreover, if $f$ is a translation, i.e. $\lambda=1$,  then we have $$|c|^2|\mu|^2 \geq \frac{(1-\rho+\sqrt{(1-\rho)^2-4\beta(f)}}{2}.$$

\end{theorem}

\medskip 

\subsection{Useful Results}\label{dstu}

Let $\mathcal L$ be the set of loxodromic elements in ${\rm U}(n,1; \F)$. 
It is well known that $\mathcal L $ is an open subset of ${\rm U}(n,1; \F)$. This fact will be crucial for our proofs. 

Let $\mathcal E$ be the set of all regular elliptic elements in ${\rm U}(n,1; \F)$.   When $\F=\C, \H$, $\mathcal E\neq \emptyset$. When $\F=\R$, note that $\mathcal E \neq \emptyset$ if and only if $n$ is even. For $n$  odd, an elliptic $f$ in ${\rm U}(n,1; \R)$ has at least two fixed points on $\partial \h^{n}$. It is known that $\mathcal E$ is an open subset of ${\rm U}(n,1; \F)$.

\medskip The following theorem will also be useful for our purpose. 

\begin{theorem}\label{cor3}\cite{cg}
 Let $G$ be a subgroup of ${\rm U}(n,1; \F)$ such that there is no point in $\overline{\fh^{n}}$ or proper totally geodesic submanifold in $\fh^{n}$ which is invariant under $G$. Then $G$ is either discrete or dense in ${\rm U}(n,1; \F)$.
\end{theorem}

\subsection{Limit set}
Let $L(G)$ be the limit set of a subgroup $G$ of ${\rm U}(n,1; \F)$. The limit set $L(G)$ is a closed $G$-invariant subset of $\partial \fh^n$. The group $G$ is elementary if $L(G)$ is finite. If $G$ is elementary, $L(G)$ consists of at most two points. If $G$ is non-elementary, then $L(G)$ is an infinite set and every non-empty, closed $G$-invariant subset of $\partial \fh^n$ contains $L(G)$. We note the following lemma, for proof see \cite[Chapter 12]{rat}.

\begin{lemma}\label{nel}
Let $a \in \partial \fh^{n}$ be fixed by a non-elliptic element of a subgroup $G$ of  ${\rm U}(n,1; \F)$, then $a$ is a limit point of $G$.
\end{lemma}

\section{Proof of \thmref{thm1} }\label{ms}
  Let $Fix(f)$ be subset of ${\overline\h}^{n+2}$ that is pointwise fixed by $f$. Let $O_f$ be the stabilizer subgroup of $Fix(f)$ in $\c$. Clearly, $O_f$ is a closed subgroup of $\c$.

If possible suppose $G$ is not discrete.  Since $G$ is Zariski-dense and assumed to be non-discrete, by \thmref{cor3},  $G$ is dense in $\c$. Let $f$ be a `test map'.  Then there exists a sequence $\{g_n\}$ of distinct  loxodromic elements such that $g_n \to f$. We may further assume that $Fix(g_n) \cap Fix(f)=\emptyset$.  Clearly, there is such a sequence $g_n'$  in $\c$.  Since $G$ is dense in $\c$,   we can choose $g_n$ sufficiently close to $g_n'$ in the open neighbourhood $\mathcal L \setminus O_{f}$.

\medskip $(1)$ Let $f$ be loxodromic.  Upto conjugacy, assume $f$ fixes $0$ and $\infty$, that is, \begin{equation} \label{lox} f= \begin{pmatrix} \lambda & 0\\ 0 & {\lambda^*}^{-1}\end{pmatrix}, \;\; |\lambda| \neq 1. \end{equation} 
Let 
\begin{equation} g_n = \begin{pmatrix}
                       a_n & b_n \\
                       c_n & d_n
                     \end{pmatrix} \end{equation}
It can be seen that $[0,\infty,g_n(0),g_n(\infty)]=-b_nc_n^*$. By \lemref{nel},  the subgroup $\langle f,g_n \rangle$ has more than two limit points, so it is  non-elementary, also discrete by hypothesis. Thus  using \thmref{watt1} and by the hypothesis, 
\begin{eqnarray*}
  & & \beta(f)(1+|b_nc_n|) \geq 1 \\
  &\Rightarrow& |b_nc_n| \geq -1 + \frac{1}{\beta(f)} >  0.
\end{eqnarray*}
But we have $b_nc_n \to 0$ as $n \to \infty$.  This leads to a contradiction.

\medskip $(2)$ Let $f$ be non-elliptic. Applying suitable conjugation, without loss of generality we may assume that one of the fixed point of $f$ be $\infty$ which leaves $f$ in the form $f = \begin{pmatrix} \lambda & \mu\\ 0 & {\lambda^*}^{-1}\end{pmatrix}$. By \lemref{nel} and hypothesis, for large $n$,  the subgroup $\langle f,g_n \rangle$ is  non-elementary and discrete. 
 Then using \thmref{watt3} we must have
$$|tr^2(fg_nf^{-1})[fg_n(\infty),fg_n^{-1}(\infty),g_n(\infty),g_n^{-1}(\infty)]| \geq  \frac{(1-\rho+\sqrt{(1-\rho)^2-4\beta(f)}}{2}\;.$$
By calculation, we see that the left hand side of the above inequality will be same as the left hand side  of the following inequality:
$$|\lambda|^{-2}|c_n|^2{|f(a_nc_n^{-1})-(a_nc_n^{-1})|.|f(-c_n^{-1}d_n)-(-c_n^{-1}d_n)|} \geq  \frac{(1-\rho+\sqrt{(1-\rho)^2-4\beta(f)}}{2},$$
i.e. $$k_n=|c_n|^2{|f(a_nc_n^{-1})-(a_nc_n^{-1})|.|f(-c_n^{-1}d_n)-(-c_n^{-1}d_n)|} \geq  \frac{|\lambda|^2 (1-\rho+\sqrt{(1-\rho)^2-4\beta(f)}}{2}.$$
Since $f$ and $g_n$ does not have a common fixed point, we must have $c_n \neq 0$. Also since $0 < \rho <1$,  hence,  $\frac{(1-\rho+\sqrt{(1-\rho)^2-4\beta(f)}}{2}$ is a positive real number. So,  $|f(a_nc_n^{-1})-(a_nc_n^{-1})|$ and $|f(-c_n^{-1}d_n)-(-c_n^{-1}d_n)|$ are non-zero. Thus for all $n$, $k_n$ is bounded above by a positive real number. But  $k_n \to 0$ as $n \to \infty$.  This is a contradiction.

\medskip $(3)$ Let $f$ be elliptic as given. Recall that in the case when $n$ is even and $f$ has no fixed points on $\partial \h^{n+2}$, we use inclusion to view $f$  as an element in ${\mathrm{ SL}}(2, C_{n+1} )$ and assume $0$, $\infty$ to be points on $\partial \h^{n+3}$, and thus 
 \begin{equation}\label{ell} f= \begin{pmatrix} \lambda & 0\\ 0 & {\lambda^*}^{-1}\end{pmatrix}, \;\; |\lambda| = 1. \end{equation} 
By hypothesis, $\langle f, g_n \rangle$ is discrete. 
We claim that $\langle f, g_n \rangle$ is non-elementary. If not, then it must keep the fixed points of $g_n$ invariant. Since $f$ does not have a common fixed point with $g_n$, it much swipes the fixed points of $g_n$. That would imply that $f$ must have a rotation angle $\pi$. But then $\beta(f)$ would be more than $4 \sin^2(\pi/10)$, which is not possible by assumption. 

Let $g_n$ be of the form \eqnref{form}. 
Since $b_nc_n \to 0$, for large $n$, 
$$\beta(f)\bigg (\frac{1}{4 \sin^2(\pi/10)}+|b_nc_n| \bigg ) < 1.$$ This is a contradiction to \thmref{watt2}. 

\medskip 
This proves the theorem.

\subsection{Proof of \thmref{thm2} }$\;\;$

\medskip 
As above, given a test map $f$ we choose a sequence of loxodromic elements $g_n$ such that $g_n \to f$ and $Fix(g_n) \cap Fix(f)=\emptyset$. Let $L_n = g_nfg_n^{-1}= \begin{pmatrix}
                           a_n & b_n \\
                           c_n & d_n
                         \end{pmatrix}$.
Note that $Fix(L_n)=g_n(Fix(f))$. 

\medskip $(1)$ Let $f$ be of the form \eqnref{lox}.  Since $g_n$ does not fix the boundary fixed points of $f$, $Fix(L_n)$ would be disjoint from $Fix(f)$. Thus $\langle f,L_n \rangle$ is non-elementary as the limit set contains $Fix(f) \cup Fix(L_n)$,  and it is discrete by hypothesis. Hence  by \thmref{watt1}, we have
  $$|b_nc_n| \geq -1+\frac{1}{\beta(f)} >0. $$
But as $L_n \to f$,  we have $b_nc_n \to 0$ as $n \to \infty$. This leads to a contradiction.

\medskip $(2)$ In this case, we follow the similar arguments as in \thmref{thm1},  and we get by \thmref{watt3} inequality that,
 $$|c_n|^2{|f(a_nc_n^{-1})-(a_nc_n^{-1})|.|f(-c_n^{-1}d_n)-(-c_n^{-1}d_n)|} \geq  \frac{|\lambda|^2(1-\rho+\sqrt{(1-\rho)^2-4\beta(f)}}{2}.$$
But since $L_n \to f$, so $c_n \to 0$, and hence $k_n \to 0$ as $n \to \infty$.
This leads to a contradiction.

\medskip $(3)$ Let $f$ be a regular elliptic. Let $f$ be of the form \eqnref{ell}. We claim that $\langle f, L_n \rangle$ is non-elementary. If not then, it either fixes a point or keeps a two point set $\{a, b\}$ on the boundary invariant. If $\langle f, L_n \rangle$  fixes a point $p$ on $\h^n$, then $f$ fixes the geodesic $l$ joining $p$ and $g_n^{-1}(p)$. Consequently $f$ fixes the boundary points of $l$. But that would imply, $f$ must preserve $g_n^{-1}(l)$. If $g_n$ does not preserve $l$, this would imply that $f$ must have another boundary fixed point or a rotation angle $\pi$, both not possible by assumption. So $g_n$ must keep $l$ invariant. This is again not possible. 

  If $\langle f, L_n \rangle$ keeps $\{a, b \}$ invariant, then $f$ keeps $g_n^{-1}(l)$ invariant, where $l$ is the geodesic joining $a$ and $b$. Thus $f$ either fixes $a$, $b$ or swipes them. If $f$ swipes them, it must have a rotation angle $\pi$ which is not possible given the value of $\beta(f)$. If $f$ fixes $a$ and $b$, then $\{a, b \}$ must be $\{0, \infty\}$. Since $L_n$ also preserves $l$,  $g_n$ must preserve $l$ joining $0$ and $\infty$. 
This is not possible because $g_n$ and $f$ do not have the same fixed points, and if $g_n$ swipes them, it must have a fixed point on $\h^n$, which is again impossible. Hence $\langle f, L_n \rangle$ must be non-elementary, and also discrete by hypothesis. Now the result follows similarly as in the proof of  \thmref{thm1}(3). 

This proves the theorem.

\section{Proof of  \thmref{thmq} }\label{qc}
 Recall that 
$$\Sp(n, 1) = \{A \in {\rm GL}(n + 1, \H) : A^*J_{2}A = J_{2}\},$$
where $$J_2=\begin{pmatrix} 0 & -1 & 0 \\ -1 & 0 & 0 \\ 0 & 0 & I_{n-1} \end{pmatrix}.$$
Equivalently, one may also use the Hermitian form given by the following matrix wherever convenient. 
$$J_1=\begin{pmatrix} -1 &0 \\ 0 & I_{n} \end{pmatrix}.$$ 

An element  $g\in{\rm Sp}(n,1)$ acts on
$\overline{\hh}^n={\bf H}_\H^n\cup\partial \hh^n$ by projective transformations.  Thus the isometry group of $\hh^n$ is given by 
${\rm PSp}(n,1)=\Sp(n,1)/\{I, -I\}.$ For a matrix (or a vector) $T$ over $\H$, let $T^{\ast}=\bar T^t$.   Let $A$ be an element in $\Sp(n,1)$. Then one can choose $A$ to be of the following form. 
\begin{equation}\label{form}
A=\begin{pmatrix} a & b & \gamma^{\ast} \\ c & d & \delta^{\ast} \\ \alpha & \beta & U \end{pmatrix},  \end{equation}
where $a, b, c, d$ are scalars, $\gamma, \delta, \alpha, \beta$ are column matrices in $\H^{n-1}$ and $U$ is an element in $M(n-1, \H)$. 
Then, it is easy to compute that 
$$A^{-1}=\begin{pmatrix} \bar d & \bar b & -\beta^{\ast} \\ \bar c & \bar a & -\alpha^{\ast} \\ -\delta & -\gamma & U^{\ast} \end{pmatrix}.$$
Let $o, \infty \in \partial \hh^n$ stand for the  vectors $(0, 1, \dots , 0)^t$ and $(1, 0, \dots, 0)^t \in \H^{n+1}$ under the projection map  respectively.

\subsection{Quaternionic hyperbolic J\o{}rgensen   inequalities} \label{shi}
For two generator subgroups of $\Sp(n,1)$ with an elliptic generator, one has the following, see \cite{ct1}, \cite{gmt2}. For elliptic elements, we use the form $J_1$ to represent $\Sp(n, 1)$. 
\begin{theorem} \label{eljo} \cite{ct1} 
Let $g$ and $h$ be elements of $\Sp(n, 1)$. Suppose that $g$ is a regular elliptic element with fixed point $0=(0,\ldots, 0)^t \in \hh^n$, i.e. $g$ is of the form 
\begin{equation}\label{ell} g=\begin{pmatrix} \lambda_1 & 0 \\ 0 & L \end{pmatrix},\end{equation} 
where $L=diag(\lambda_2, \ldots, \lambda_{n+1})$. Let 
$$ h = (a_{i,j})_{i, j = 1, \dots , n+1} = \begin{pmatrix}    a_{1, 1} & \beta \\ \alpha & A \end{pmatrix},$$
be an arbitrary element in $\Sp(n, 1)$, where $a_{1, 1}$ is a scalar, $\alpha, \beta$ column vectors and $A \in M(n, \H)$.  If 
$$|a_{1, 1}| \delta(g)<1.$$
then the group $\langle g, h \rangle$ generated by $g$ and $h$ is either elementary or non-discrete.
\end{theorem}

For representing  parabolic and loxodromic elements, we shall use the Hermitian form $J_2$. 
 In  \cite[Appendix]{hp}, Hersonsky and Paulin proved a version of Shimizu's lemma  for subgroups in $\SU(n,1)$. The following quaternionic version of \cite[Proposition A.1]{hp} is a straight-forward adaption of the proof of Hersonsky and Paulin. 
\begin{theorem}\label{sht}
Suppose $T_{s, \zeta}$ be an Heisenberg translation in $\Sp(n,1)$ of the form \eqnref{ht1},  and $A$ be an element in $\Sp(n,1)$ of the form \eqnref{form}.  Set
\begin{equation}t={\rm Sup}\{|b|, |\beta|, |\gamma|, |U-I| \}, \ M=|s|+2|\zeta|. \end{equation}
If 
\begin{equation}M t +2|\zeta| < 1,\end{equation}
then the group generated by $A$ and $T_{s, \zeta}$ is either non-discrete or fixes $o$. 
\end{theorem}

For two generator subgroups with a loxodromic element, we have the following version of the J\o{}rgensen inequality from the work of Cao and Parker \cite{cp1}. Up to conjugacy, a loxodrmic element has fixed points $o$ and $\infty$, and it is conjugate to a matrix of the form \eqnref{elpl}. 
\begin{theorem} {\rm (Cao and Parker)  \cite{cp1}}\label{cpt}  Let $h \in {\rm Sp}(n, 1)$ be given by (\ref{form}).
Let $g$ be  a loxodromic element in $\Sp(n,1)$ with fixed points $o,  \infty  \in \partial \hh^n$, i.e. of the form \eqnref{elpl}. Let  $M_g<1$. If $\langle g, h \rangle$ is non-elementary and discrete, then 
\begin{equation}\label{eql} |ad|^{\frac{1}{2}}|bc|^{\frac{1}{2}} \geq \frac{1 - M_g}{M_g^2}. \end{equation} 
\end{theorem}

\subsection{Proof of \thmref{thmq} } \label{pfthd}

If possible suppose $G$ is not discrete. Then $G$ must be dense in $\Sp(n,1)$ by \thmref{cor3}. 
Note that the set $\mathcal L$ of loxodromic elements in $\Sp(n,1)$ forms an open subset of $\Sp(n,1)$. Let  $Fix(f)$ denote the fixed point set of $f$ on $\partial \hh^n$. Let $F_f$ be the subgroup of $\Sp(n,1)$ that stabilizes $Fix(f)$. The subgroup $F_f$ is closed in $\Sp(n,1)$. Hence $\mathcal L - F_f$ is still an open subset in $\Sp(n,1)$. 

\medskip (1) Let $f$ be loxodromic. Up to conjugacy, assume that $f$ is of the form \eqnref{elpl}. 
Since $f \in \bar G$,  using similar arguments as in the proof of \thmref{thm1}, there exists a sequence $\{h_n\}$ of loxodromic elements in $(\mathcal L - F_f) \cap G$ such that $h_n \rightarrow f$.  Thus, $h_n, f$ do not have a common fixed point, and $\langle h_n, f \rangle$ is non-elementary for each $n$. Let
$$ h_n = \begin{pmatrix}
     a_n & b_n & \gamma_n^*\\
c_n  & d_n & \nu_n^*\\
    \alpha_n &  \beta_n & U_n\\
     \end{pmatrix},$$
where $a, b, c, d$ are scalars, $\gamma, \delta, \alpha, \beta$ are column matrices in $\H^{n-1}$ and $U$ is an element in $M(n-1, \H)$. By \thmref{cpt}, 
    $$|a_nd_n|^{\frac{1}{2}} |b_nc_n|^{\frac{1}{2}} > \frac{1 - M_f}{M_f^2}.$$
    But $b_nc_n \rightarrow 0$ as $n \to \infty $, hence 
    $$\frac{1 - M_f}{M_f^2} < 0, $$
    which is a contradiction since $M_f <1$. 

\medskip (2) Let $f$ be a Heisenberg translation. Without loss of generality assume it is of the form \eqnref{ht1}. Since, $f \in \overline{G}$, there exist a sequence of loxodromic elements $\{h_n\} \in (\mathcal L - F_f)\cap G$ such that 
$$ h_n \rightarrow f.$$
Since, $f$ and $h_n$ have distinct fixed points, hence  $\langle f, h_n \rangle$ is discrete and non-elementary. By \thmref{sht}, 
$$Mt_n + 2 |\zeta| \geq 1.$$
But $t_n \to 0$ as $n \to \infty$. Thus for large $n$, $|\zeta|\geq \frac{1}{2}$. 
This is a contradiction as $|\zeta| < \frac{1}{2}$ is given.  

\medskip 
(3) Let $f$ be a regular elliptic. We can assume that $f$ is of the form \eqnref{ell} with fixed point $0$,  up to conjugacy. Since, $G$ is dense in $\Sp(n,1)$, there is  a sequence of loxodromic element $\{h_m\}$ in $\mathcal L \cap G$ such that $h_m \to I$.  Let 
$$h_m = (a_{i,j}^{(m)}) = {\begin{pmatrix}    a_{1, 1}^{(m)} & \beta^{(m)}\\\alpha^{(m)} & A^{(m)} \end{pmatrix}}.$$
The group $\langle f, h_m\rangle$ must be non-elementary. For, if not, clearly $\langle f, h_m \rangle$ can not fix a point on $\overline \hh^n$ as that will contradict either regularity of $f$ or loxodromic nature of $h_m$.  If it keeps two points $x$ and $y$ on $\partial \hh^n$ invariant without fixing them, then $f$ must swipes $x$ and $y$, and hence  $f^2$ fixes $x$ $y$, and $0$.  Thus $f^2$ must have a repeated eigenvalue $\lambda$, see \cite[ Proposition 2.4]{cg}. This implies, $g$ would have a repeated eigenvalue $\lambda^{1/2}$, which is a contradiction to the regularity.  By our assumption $\langle f, h_m \rangle$ is also discrete for each $m$. Hence by \thmref{eljo}, 
$$|a_{1, 1}^{(m)}| ~\delta(g) \geq 1. $$
But $a_{1, 1}^{(m)} \to 1$ and $\delta(g) <1$.  This is a contradiction. 

This proves the theorem.

\begin{remark}
The results in this paper show that in order to determine discreteness of a Zariski-dense subgroup $G$ of ${\rm U}(n, 1;\F)$,  it is enough to check discreteness of the two generator subgroups of $G$ obtained by adjoining the loxodromic elements of $G$  to a `test map'  in ${\rm U}(n, 1;\F)$.  Let $\mathcal E$ denote the set of regular elliptic elements of ${\rm U}(n, 1;\F)$. The set $\mathcal E$ is also an non-empty open subset of ${\rm U}(n, 1;\F)$, provided $n$ is even when $\F=\R$. Thus, if we replace the loxodromic elements $g$ by regular elliptic elements, then versions of \thmref{thmq} and \corref{cord} hold true for all $n$, and, \thmref{thm1}  goes through for all even $n$. 
\end{remark} 

\end{document}